\documentclass[11pt]{amsart}

\def\proof{{\sc Proof. }}
\def\endproof{\hfill\hskip -1cm \rule{.5em}{.5em}}
\def\R{\mathbf{R}}

\def\lie{\mathcal{L}}
\def\nablaM{\stackrel{M}{\nabla}}
\def\nablaN{\stackrel{N}{\nabla}}

\newcommand{\gata}{\hfill\hskip -1cm \rule{.5em}{.5em}}

\newtheorem{contor}{1}[section]
\newtheorem{defi}[contor]{Definition}
\newtheorem{prop}[contor]{Proposition}
\newtheorem{cor}[contor]{Corollary}
\newtheorem{theo}[contor]{Theorem}
\newtheorem{rem}[contor]{Remark}
\newtheorem{ex}[contor]{Example}

\title[Gray curvature identities]
    {\bf  Gray Curvature Identities for Almost Contact Metric Manifolds}
\thanks{This work was supported by Grant CEEX ET n. 5883/2006-2008 ANCS Romania}
\author{Raluca Mocanu {\tiny and} Marian Ioan Munteanu}
\date{\today}

\begin{document}

\maketitle
\begin{abstract}
The aim of this research is the study of Gray curvature identities,
introduced by Alfred Gray in \cite{kn:Gra76} for the class of almost
hermitian manifolds. As known till now, there is no equivalent for
the class of almost contact manifolds.
For this purpose we use the Boohby-Wang fibration and the warped manifolds
construction in order to establish which identities could be satisfied by
an almost contact manifold.
An almost hermitian manifold which satisfies one of the three Gray
identities has rich topological and geometric properties.\\[1mm]
{\bf Keywords and Phrases:} almost Hermitian manifolds, almost contact metric manifolds,
curvature identities, Boothby Wang fibration, cone metric, cosymplectic manifolds,
Sasakian manifolds, generalized Heisenberg group.\\[1mm]
{\bf Mathematics Subject Classification (2000):} 53C15, 53C25, 53C55, 53B35, 53D15.
\end{abstract}
\medskip

\section{Introduction}
In their paper \cite{kn:BHR90}, the authors defined $K_{i\varphi}$-curvature
identities $(i=1,2,3)$ for an almost contact metric manifold $(M,\varphi,\xi,\eta,g)$
by using the usual Hermitian structure on $M\times\R$ (the product manifold).
It is known that both cosymplectic and Sasakian manifolds are natural
odd-dimension versions for Kaehlerian manifolds. Gray proved in \cite{kn:Gra76}
that Kaehlerian manifolds satisfy $K_i$, $i=1,2,3$ (curvature identities for almost
Hermitian manifolds). In the same spirit, in \cite{kn:BHR90} it is shown that
cosymplectic manifolds satisfy $K_{i\varphi}$-identities. We asked what happens
with Sasakian manifolds? Recall that a Riemannian manifold $(M,g)$ is Sasakian
if the holonomy group of the metric cone on $M$:
$(C(M)=\R_+\times M,\widetilde g=dt^2+t^2g)$ reduces to a subgroup of
$U\left(\frac{m+1}2\right)$, i.e. $(C(M),\widetilde g)$ is Kaehlerian.
(Here $m=\dim M$.) Inspired from this definition and from \cite{kn:BHR90}
we will give another approach of Gray curvature identities for almost contact
metric manifolds.

\subsection{Gray curvature identities}
An almost Hermitian manifold $(M,J,g)$ is said to satisfy the Gray curvature identities
$(K1)$, $(K2)$ and respectively $(K3)$,
if his Riemann Christoffel curvature tensor verifies \\[2mm]
$(K1)\qquad
      R(X,Y,Z,W) = R(X,Y,JZ,JW)$\\[2mm]
$(K2)\qquad
     R(X,Y,Z,W) = R(JX,JY,Z,W)+ R(JX,Y,JZ,W) + R(JX,Y,Z,JW)$\\[2mm]
 $(K3)\qquad    R(X,Y,Z,W) = R(JX,JY,JZ,JW)$\\[2mm]
for all vector fields $X,Y,Z,W$ on $\chi(M)$.
Throughout of this paper, the curvature tensor
is defined by $R_{XY}Z = \nabla _X \nabla _Y Z - \nabla _Y \nabla _X Z - \nabla _{[X,Y]} Z$,
for all $X,Y,Z\in\chi(M)$ while the Riemann Christoffel curvature tensor is given by
$R(X,Y,Z,W) = -g(R_{XY}Z,W)$.

\section{Warped product manifolds} {\em Singly warped products} or simply {\em warped products} were first
defined by Bishop \& O'Neill in \cite{kn:BON} in order to construct Riemannian
manifolds with negative sectional curvature. Let $(B,g_B)$ and $(F,g_F)$ be Riemannian
manifolds and let $b:B\longrightarrow (0,\infty)$ be a smooth function.
The warped product $\widetilde M=B\times_bF$ is the product manifold $B\times F$
endowed with the metric $\widetilde g=g_B\oplus b^2g_F$. More
precisely, if $\pi : B\times F\longrightarrow B$ and
$\tau : B\times F\longrightarrow F$ are natural projections, the metric
$g$ is defined by
\begin{equation}
\label{eq3}
\widetilde g=\pi^*g_B+(b\circ\pi)^2\tau^*g_F.
\end{equation}
The function $b$ is called {\em warping function}.
If $b\equiv1$, then we have a product manifold.

\medskip
If $X,Y$ are tangent to $B$ and $Z,W$ tangent to $F$, then the Levi Civita
connection $\widetilde\nabla$ of $\widetilde M$ is given by
\begin{equation}
\label{eq4}
\left\{
\begin{array}{l}
\widetilde \nabla_XY=\nabla_X^BY,\quad
\widetilde\nabla_XZ=X(\ln b)Z\\[1mm]
\widetilde\nabla_ZW=\nabla_Z^FW-b^2\ g_F(Z,W)\nabla^B(\ln b)
\end{array}\right.
\end{equation}
where $\nabla^B$ and $\nabla^F$ are the Levi Civita connections on $B$,
respectively on $F$, and $\nabla^B(\ln b)$ is the gradient of $\ln b$
with respect to the metric $g_B$.

\medskip
Let $(M, \varphi, \xi, \eta, g)$ be an  almost contact metric manifold.
Consider the warped product manifold $\widetilde M = \R_+ \times_t M$, where $t$ is
the global coordinate of $\R_+$, i.e. the metric $\widetilde g$ of $\widetilde M$
is defined by
\begin{equation}
\label{k1}
%\usecounter{contor}
\widetilde g = dt^2 + t^2g.
\end{equation}
Define an endomorphism on $\chi(\widetilde M)$ by
\begin{equation}
\label{k2}
%\usecounter{contor}
J \partial_t=-\frac{1}{t}\ \xi
\:\:\: JX=\varphi X +t\eta(X) \partial_t,   \quad \forall X \in \ \chi (M)
\end{equation}
where $\partial_t=\frac{d}{dt}$.
For $\widetilde X=(a,X) \in \chi(\widetilde M)$, $a \in C^\infty (\R_+), X \in \chi(M)$ we have
\begin{equation}
\label{k3}
%\usecounter{contor}
J \widetilde X=J(a,X)=\big(t \eta(X), \varphi X - \frac {a}{t}\ \xi\big)
\end{equation}
The proofs of the following propositions are straightforward.
\begin{prop}
$J$ is an almost complex structure
compatible with the metric $\widetilde g$.
\end{prop}
\begin{prop}
The Levi-Civita connection $\widetilde \nabla$ of $\widetilde g$ is given by:
\begin{equation}
%\usecounter{contor}
\left\{
  \begin{array}{l}
         \widetilde \nabla_{\partial_t} \partial_t =0,\quad
         \widetilde \nabla_X \partial_t = \widetilde \nabla_{\partial_t} X = \frac{1}{t} X  \\[2mm]
         \widetilde \nabla_X Y = \nabla_X Y -tg(X,Y) \partial_t,      \quad X,Y \in \chi(M)
 \end{array}
\right.
\end{equation}
\end{prop}

\begin{prop}
The covariant derivative of J is given by:
\begin{equation}
%\usecounter{contor}
   \left\{
     \begin{array}{l}
(\widetilde \nabla_{\partial_t} J) \partial_t =(0,0),  \quad (\widetilde \nabla_{\partial_t}) X =(0,0)\\[2mm]
(\widetilde \nabla_X J) \partial_t = (0, -\frac{1}{t} (\nabla_X \xi + \varphi X)) \\[2mm]
(\widetilde \nabla_X J) Y = (t((\nabla_X \eta)(Y) -g(X, \varphi Y)), (\nabla_X \varphi) Y-g(X,Y) \xi + \eta(Y) X)
\end{array}\right.
\end{equation}
\end{prop}

\begin{cor}
J is parallel if and only if
\begin{equation}
%\usecounter{contor}
\left\{\begin{array}{l}
(\widetilde \nabla_X \varphi) Y =g(X,Y)\xi - \eta(Y)X,  \quad
(\nabla_X \eta)(Y) = g(X, \varphi Y) \\[2mm]
\nabla_X \xi = - \varphi X,   \:\:\:\:X, Y \in \chi(M)
\end{array}\right.
\end{equation}
i.e. $(\widetilde M, J, \widetilde g)$ is Kaehler if and only if $(M,\varphi, \xi, \eta, g)$ is Sasakian.
\end{cor}
\begin{prop}
For the curvature of the manifold $\widetilde M$ we have
\begin{equation}
%\usecounter{contor}
\left\{\begin{array}[c]{lcr}
 \widetilde R(\partial_t, X) \partial_t =0,\
 \widetilde R(X,Y) \partial_t =0, \
 \widetilde R(\partial_t,X) Y =0 \\[2mm]
 \widetilde R(X,Y) Z = R(X,Y)Z - g(Y,Z)X + g(X,Z)Y
\end{array}\right.
\end{equation}
where $\widetilde R$  (respectively $R$) are the curvature tensors for $\widetilde g$
(respectively for $g$).
\end{prop}
\begin{prop}
Moreover, the following relations hold:
\begin{equation}
%\usecounter{contor}
\left\{\begin{array}[c]{lcr}
 \widetilde R(\partial_t, X) (J \partial_t )=0,\  \widetilde R(\partial_t, X) (J Y )=0\\[2mm]
 \widetilde R(X,Y) (J \partial_t )= -\frac{1}{t} [R(X,Y) \xi - \eta(Y) X + \eta(X)Y]\\[2mm]
 \widetilde R(X,Y) (JZ) =R(X,Y) (\varphi Z)-g(Y, \varphi Z) X + g(X, \varphi Z) Y
\end{array}\right.
\end{equation}
\end{prop}

In the following we compute expressions of the form
$\widetilde g (\widetilde R(A,B) (JC), JD)$.
The useful expressions are obtained in the following cases:\\[1mm]

{\bf 1.}
$ \widetilde g (\widetilde R(X,Y) (J \partial_t ), JW)=
     -t\left[ g(R(X,Y) \xi, \varphi W) - \eta(Y) g(X, \varphi W) + \eta(X) g(Y, \varphi W)\right]
     $

{\bf 2.}
$  \widetilde g (\widetilde R(X,Y) (J Z),J \partial_t)=
     -t\left[ \eta(R(X,Y)(\varphi Z)) - \eta(X) g(Y, \varphi Z) + \eta(Y) g(X, \varphi Z)\right]
$

{\bf 3.}
$
 \widetilde g( \widetilde R(X,Y) (JZ), JW)=
        t^2\left[ g(R(X,Y) \varphi Z, \varphi W)- g(Y, \varphi Z) g(X, \varphi W) +\right.$

\qquad\qquad\qquad\qquad\qquad\qquad\qquad\qquad
  $\left. +g(X, \varphi Z)g(Y, \varphi W)\right]$.

\begin{theo} $\widetilde M$ is $(K1)$ if and only if
\begin{equation}\label{k4}
%\usecounter{contor}
\begin{array}{c}
 R(X,Y,Z,W) = R(X,Y, \varphi Z, \varphi W) -g(X, \varphi Z)g(Y, \varphi W)+\\[2mm]
\qquad\qquad\qquad\qquad +g(Y, \varphi Z)g(X, \varphi W)-g(Y,Z)g(X,W)+g(X,Z)g(Y,W).
\end{array}
\end{equation}
\end{theo}
\proof $\widetilde M$ is $(K1)$ if and only if
$\widetilde g(\widetilde R(A,B)(JC),JD)=\widetilde g(\widetilde R(A,B)C,D)$\\[1mm]
for all $A,B,C,D \in \chi(\widetilde M)$
\medskip

${\mathbf{1}.}\
       \widetilde g(\widetilde R(X,Y) (J \partial_t),JW)=\widetilde g(\widetilde R(X,Y) \partial_t,W)
       $

\qquad $\Longrightarrow$
$-t[g(R(X,Y) \xi, \varphi W) - \eta(Y) g(X, \varphi W) + \eta (X) g(Y, \varphi W)]=0.$

\qquad
$\Longrightarrow$ $g(\varphi(R(X,Y) \xi - \eta (Y)X+\eta (X)Y,W)=0$, for every W

\qquad $\Longrightarrow$
$R(X,Y)\xi - \eta(Y)X+ \eta(X)Y \in \ker \varphi$.

Thus, we have obtained
\begin{equation}
\label{k5}
%\usecounter{contor}
R(X,Y) \xi = \eta (Y)X-\eta (X)Y\:\:\: modulo\:\: \xi.
\end{equation}

${\mathbf 2.}\
 \widetilde g(\widetilde R(X,Y)(JZ),J \partial_t)= \widetilde g(\widetilde R(X,Y) Z, \partial_t)$

\qquad $\Longrightarrow$
$\widetilde g\left(R(X,Y)(\varphi Z) - g(Y, \varphi Z) X + g(X, \varphi Z)Y, -\frac{1}{t}\ \xi \right)=$

\qquad \qquad\qquad\qquad\qquad
$ =\widetilde g \left(R(X,Y)Z - g(Y,Z)X+g(X,Z)Y, \partial_t\right)$

\qquad $\Longrightarrow$
$
-\frac{1}{t}\ t^2 \eta(R(X,Y) (\varphi Z) -g(Y, \varphi Z)X + g(X,\varphi Z) Y)=0$.

\medskip
Thus,
\begin{equation}\label{k6}
%\usecounter{contor}
R(X,Y)(\varphi Z)-g(Y, \varphi Z)X+g(X, \varphi Z)Y \in Ker\: \eta
\end{equation}

${\mathbf 3.}\
\widetilde g(\widetilde R(X,Y) (JZ), JW)= \widetilde g(\widetilde R(X,Y)Z,W)$

\qquad $\Longrightarrow$
$\widetilde g(R(X,Y)(\varphi Z)-g(Y, \varphi Z)X+g(X,\varphi Z)Y, \varphi W+
  t\eta (W) \partial_t)=$

\qquad\qquad\qquad\qquad
  $=\widetilde g(R(X,Y)Z-g(Y,Z)X+g(X,Z)Y,W)$.

After simplification by $t^2$ we obtain
$$g(R(X,Y)(\varphi Z), \varphi W) -g(Y, \varphi Z)g(X, \varphi W) +g(X, \varphi Z) g(Y, \varphi W)=$$
$$=g(R(X,Y)Z,W)-g(Y,Z)g(X,W)+g(X,Z)g(Y,W).$$
It follows
\begin{equation}
\label{k7}
%\usecounter{contor}
\begin{array}{c}
R(\varphi W, \varphi Z,X,Y) - R(W,Z,X,Y)=g(Y, \varphi Z)g(X, \varphi W)-\quad\\
\qquad\quad  -g(X,\varphi Z)g(Y, \varphi W)+g(X,Z)g(Y,W)-g(Y,Z)g(X,W)
\end{array}
\end{equation}
\endproof

\begin{rem} \rm We immediately obtain
(\ref{k5}) $\longrightarrow$  (\ref{k6}) and (\ref{k7}) $\longrightarrow$ (\ref{k5}).
\end{rem}

Return to the formula (\ref{k7}).
We interchange $(\varphi W, \varphi Z)  \longleftrightarrow (X,Y)$,
$(W,Z) \longleftrightarrow (X,Y)$
and then $Z\longleftrightarrow W$. One gets
$$ \bf (*)\ R(X,Y,\varphi Z,\varphi W) -R(X,Y,Z,W) = g(Y, \varphi W)g(X,\varphi Z) -g(X, \varphi W)g(Y, \varphi Z) +$$
$$\bf +g(X,W)g(Y,Z) -g(Y,W)g(X,Z)$$
for all $X,Y,Z,W$ in $\chi(M)$.
\medskip

As consequences we have

\quad
$R(\xi,Y,\xi,W)=g(Y,W)$

\quad
$R(\xi,Y,Z,W)=R(\xi,Y,\varphi Z,\varphi W)=0$

\quad
$R(X,Y,Z,W)-g(Y,W)g(X,Z)+g(X,W)g(Y,Z) =$

\qquad\qquad\qquad
$=R(X,Y,\varphi Z,\varphi W)-g(Y,\varphi W)g(X, \varphi Z)+g(X, \varphi W)g(Y, \varphi Z)$,

where $X,Y,Z$ and $W$ are orthogonal to $\xi$.

\vspace{2mm}
\begin{defi}
We say that an almost contact metric manifold satisfies {\bf (G1)}-identity
if its curvature tensor verifies $\mathbf{(*)}$.
\end{defi}
\begin{prop} The curvature tensor of a Sasakian manifold satisfies {\bf (G1)}
{\rm (see also Lemma 7.1 in \cite{kn:Bla02})}.
\end{prop}
\begin{prop}
\label{prop:cont_sas}
Any contact manifold satisfying {\bf(G1)} is Sasakian.
\end{prop}
\proof
It is known (e.g. Proposition 7.6 from \cite{kn:Bla02}) that a contact manifold
is Sasakian if and only if $R(X,Y)\xi=\eta(Y)X-\eta(X)Y$, for all $X$ and $Y$.
\gata

% ===============================================
% ========================== K2 =================
% ===============================================

\medskip

Return to the cone manifold $\widetilde M$. We give
\begin{theo} $\widetilde M$ is $(K2)$ if and only if
\begin{equation}
%\usecounter{contor}
R(X,Y,Z,W)= R(\varphi X, Y,Z, \varphi W)+R(X, \varphi Y,Z, \varphi W)+R(X,Y, \varphi Z,\varphi W)$$
$$+g(X,Z)\eta (W) \eta(Y) - g(Z,Y) \eta(X) \eta(W)
\end{equation}
\end{theo}
\proof
$\widetilde M$ is $(K2)$ if and only if
$$\widetilde R(A,B,C,D)= \widetilde R(JA, B,C,JD)+\widetilde R(A, JB,C, JD)+\widetilde R(A,B, JC,JD)$$

Three cases are essential:\\
1) $A = \partial_t$, $B=Y$, $C=\partial_t$, $D=W$ which is equivalent to $0=0$.\\
2) $A=\partial_t$, $B=Y$, $C=Z$, $D = W$\\[1mm]

One has:

\qquad $\widetilde R(J\partial_t,Y,Z,JW) = -\frac{1}{t} \widetilde R(\xi, Y,Z, \varphi W)$

\qquad $\widetilde R(\partial_t, JY, Z,JW)=0$

\qquad $\widetilde R(\partial_t, Y, JZ,JW)=0$

It follows that the right side is equal to:
$$-tg(\xi, R(Z, \varphi W) Y - g(\varphi W, Y) Z + g(Z, Y) \varphi W)$$
Since the left side vanishes, in this case we obtain
\begin{equation}
\label{k8}
%\usecounter{contor}
R(\xi,Y,Z,\varphi W) =\eta(Z)g(\varphi W, Y)   {\rm \ for \ every \ } Y, Z, W  \in \chi(M)
\end{equation}

3) $A=X$, $B=Y$, $C=Z$, $D=W$. One has
 $$\widetilde R(JX,Y,Z,JW) = \widetilde R(\varphi X, Y, Z, \varphi W)$$
    $$\widetilde R(X,JY,Z,JW) = \widetilde R( X, \varphi Y, Z, \varphi W)$$
    $$\widetilde R(X,Y,JZ,JW) = \widetilde R( X, Y, \varphi Z, \varphi W)$$
It follows that the right side is equal to
$$ t^2[R(\varphi X,Y,Z,\varphi W) + R(X,\varphi Y, Z, \varphi W) +R(X,Y, \varphi Z, \varphi W)]+$$
$$+t^2[-g(\varphi W,\varphi Y)g(X,Z) +g(Z,Y)g(\varphi X,\varphi W)]$$
while the left side equals to:
$$t^2 R(X,Y,Z,W) + t^2 [-g(W,Y)g(X,Z) +g(Z,Y) g(X,W)]$$
We get
$$(**)
%\label{k9}
%\usecounter{contor}
\begin{array}{c}{\mathbf{
R(X,Y,Z, W) = R(\varphi X, Y, Z,\varphi W)+R(X, \varphi Y, Z, \varphi W)+R(X,Y,\varphi Z,\varphi W)+}}\\
{\mathbf{\qquad \qquad +g(X,Z)\eta(W)\eta(Y) -g(Z,Y)\eta(X)\eta(W)}}.
\end{array}
$$

It can be proved that previous relation implies \ref {k8}. Hence the statement.

\gata

As consequences one has

\qquad $R(\xi,Y,\xi,W)=g(Y,W)$

\qquad $R(\xi,Y,Z,W)=0$

\qquad $R(X,Y,Z,W)=R(\varphi X,Y,Z, \varphi W)+R(X, \varphi Y, Z, \varphi W)+R(X,Y,\varphi Z, \varphi W)$

for all $X,Y,Z,W$ orthogonal to $\xi$.

\begin{defi}
We say that an almost contact metric manifold satisfies {\bf (G2)}-identity
if its curvature tensor verifies $\mathbf{(**)}$.
\end{defi}

\begin{theo} The manifold $\widetilde M$ is $(K3)$ if and only if
\begin{equation}
\label{k10}
%\usecounter{contor}
\begin{array}{c}
R(X,Y,Z, W) = R(\varphi X, \varphi Y, \varphi Z,\varphi W)+g(X,Z)\eta(W)\eta(Y)-\qquad\qquad \\
\qquad -g(Z,Y)\eta(X)\eta(W)+g(Y,W) \eta(X)\eta(Z) -g(X,W)\eta(Y)\eta(Z)
\end{array}
\end{equation}
for all $X,Y,Z,W \in \chi(M) $.
\end{theo}
\proof $\widetilde M$ is $(K3)$ iff $\widetilde R(A,B,C,D) = \widetilde R(JA,JB,JC,JD)$
for all $A,B,C,D \in \chi(\widetilde M).$

The essential cases are:\\[1mm]
1) $A = \partial_t$, $B=Y$, $C=\partial_t$, $D=W$.

The left member vanishes and the right member is equal to
$R(\xi, \varphi Y, \xi, \varphi W) -g(\varphi W, \varphi Y)$.
We get
\begin{equation}\label{k11}
%\usecounter{contor}
R(\xi, \varphi Y, \xi, \varphi W) = g(\varphi W, \varphi Y)
\end{equation}

2) $A=\partial_t$, $B=Y$, $C=Z$, $D=W$.

 The left member vanishes and the right member is equal to
$R(\xi, \varphi Y, \varphi Z, \varphi W)$.
We get
\begin{equation}
\label{k12}
%\usecounter{contor}
R(\xi, \varphi Y, \varphi Z, \varphi W) = 0
\end{equation}

3) $A=X$, $B=Y$, $C=Z$, $D=W$.

 The left member is equal to
$$t^2[R(X,Y,Z,W)-g(Z,X)g(W,Y)+g(Y,Z)g(X,W)]$$
and the right member is equal to
$$
  t^2[R(\varphi X, \varphi Y, \varphi Z, \varphi W) -g(\varphi W, \varphi Y)
     g(\varphi X, \varphi Z)+g(\varphi Y, \varphi Z)g(\varphi X, \varphi W)]
$$
Hence (\ref{k10}) is proved.
Remark that (\ref {k10}) implies both (\ref {k11}) and (\ref{k12}).

\gata

As consequences we have

\qquad $R(\xi,Y,\xi,W)=g(Y,W)$

\qquad $R(\xi,Y,Z,W)=0$

\qquad $R(X,Y,Z,W) = R(\varphi X, \varphi Y,\varphi Z,\varphi W)$

for all $X,Y,Z,W \in \chi(M)$ orthogonal to $\xi$.

\begin{defi}
We say that an almost contact metric manifold satisfies $\bf G3$-identity
if its curvature verifies relation {\rm (\ref {k10}).}
\end{defi}

\section{The Boothby Wang fibration} Let $M$ a
($2n+1$)-dimensional smooth manifold. A {\em contact form} on $M$ is
a $1-$form $\eta$ satisfying
$$
    \eta\wedge (d\eta)^n\neq 0.
$$
We say that $\eta$ endows on $M$ a {\em contact structure}. It is
clear that $\eta$ induces an orientation on $M$ and hence there
exists a global non vanishing vector field $\xi$ on $M$ such that
$\eta(\xi)=1$. If $\xi$ is {\em regular} in the sense of Palais (see
\cite{kn:Pal57}), then the contact structure (and also $M$) is
called {\em regular}. If moreover $M$ is compact, one can consider
the space of all orbits of $\xi$, i.e. $N=M_{/_\xi}$ obtaining a
smooth manifold. We have {\bf Theorem A} (\cite{kn:BW58}). {\em Let
$(M,\eta)$ be a compact, regular, contact manifold. Then $M$ is a
principal circle bundle over $N$ and $\eta$ is a connection form of
this bundle. The curvature form $\Theta$ of $\eta$ defines a
symplectic form on $N$.} This fibration $S^1\longrightarrow
M\stackrel{\pi}{\longrightarrow} N$ is called
the Boothby-Wang fibration. \\[2mm]
Let $\Omega $ the symplectic 2-form of N, we denote by G the
associated metric, i.e.
$\Omega (X,Y)=G(X,JY)$ with $J$ the almost complex structure.\\[2mm]
In the following, we denote by by $X^\uparrow$ the lift of a vector
field $X\in \chi (N)$. $X^\uparrow$ is a horizontal vector field of
$M$. On $M$ a $(1,1)$ tensor field $\varphi$ can be defined, namely
\begin{equation}
    \varphi X^\uparrow = (JX)^\uparrow\quad , \quad \varphi\xi=0.
\end{equation}
We can easily see that
$$
\varphi ^2=-I+\eta \otimes \xi
$$
In this way, $(\varphi ,\xi ,\eta )$ becomes an almost contact
structure. The metric $G$ can be lifted and hence one defines $g$ on
$M$ as follows:
\begin{equation}
g=\pi ^*G+\eta \otimes \eta
\end{equation}
The metric $g$ is compatible with the contact structure and $\xi=\eta^\#$.\\[2mm]
Without loss of the generality one can suppose $d\eta =\pi ^*\Omega
$ and thus we have
$$
g(X^\uparrow,\varphi Y^\uparrow)=
          G(X,JY)\circ \pi =\Omega (X,Y)\circ \pi =\pi ^*\Omega (X^\uparrow,Y^\uparrow)=d\eta (X^\uparrow,Y^\uparrow)
$$
In this way,
$(\varphi ,\xi ,\eta ,g)$ becomes a contact metric structure on $M$.\\[2mm]
If the symplectic structure of $N$ derives from a Kaehlerian
structure $(J,G)$, the obtained structure on $M$ is Sasakian (i.e.
contact and normal manifold). See e.g.\cite{kn:Bla02}. But
generally, a symplectic structure need not come from a Kaehlerian
one. Yet, one can always find an almost Kaehlerian structure
inducing it. In this case, the contact structure on the total space
of a Boothby Wang fibration is $K$-contact, i.e. the vector field
$\xi$ is Killing, namely $\lie_\xi g=0$. It easily follows that the
integral curves of $\xi$ are geodesics.

It is easy to prove the relation
\begin{equation}
[X^\uparrow,Y^\uparrow]=[X,Y]^\uparrow-2G(X,JY)\xi
\end{equation}
for all $X,Y\in\chi(N)$.\\[2mm]
Denote by $\nablaM$ and $\nablaN$ the Levi Civita connections on $M$
and $N$, respectively. We immediately have:
$$
g(\nablaM_{X^\uparrow} Y^\uparrow,Z^\uparrow)\circ\pi= G(\nablaN_XY
,Z)
$$
for any $X,Y,Z \in \chi (N)$. For the vertical part %(i.e. along $\xi$)
we shall compute $\eta(\nablaM_{X^\uparrow} Y^\uparrow)$:
$$
\begin{array}{rl}
2g(\nabla^M_{X^\uparrow}  Y^\uparrow,\xi ) & = X^\uparrow
g(Y^\uparrow,\xi )+
   Y^\uparrow g(X^\uparrow,\xi )-\xi g(X^\uparrow, Y^\uparrow)+
   g([X^\uparrow,Y^\uparrow],\xi )+\\[1mm]
   & \qquad\qquad + g([\xi ,X^\uparrow],Y^\uparrow)+ g(X^\uparrow,[\xi ,Y^\uparrow])=\\[2mm]
   & = \eta ([X^\uparrow,Y^\uparrow]) - \left(\lie_\xi g\right)(X^\uparrow,Y^\uparrow)\\[2mm]
   & = -2d\eta(X^\uparrow,Y^\uparrow).
\end{array}
$$
We obtain that
$$
\eta(\nablaM_{X^\uparrow} Y^\uparrow)\circ\pi=-G(X,JY).
$$
In the following, we will ignore $\pi$, due to the isomorphism
between the horizontal distribution of $T(M)$ and $T(N)$. Hence
\begin{equation}
\nablaM_{X^\uparrow} Y^\uparrow=(\nablaN_X Y)^\uparrow -G(X,JY)\xi.
\end{equation}
In the same way, one can show
\begin{equation}
\nablaM_{X^\uparrow}\xi=-\varphi{X^\uparrow}.
\end{equation}

Denote by $R^M$ and $R^N$ the curvature tensors of $M$ and $N$,
respectively.

Then
$$
\begin{array}{c}
R^M(X^\uparrow,Y^\uparrow)Z^\uparrow
=\left(R^N(X,Y)Z\right)^\uparrow+
     g(Y^\uparrow,\varphi Z^\uparrow) \varphi X^\uparrow-
     g(X^\uparrow,\varphi Z^\uparrow) \varphi Y^\uparrow-\qquad\\[1mm]
\qquad     -2g(x^\uparrow,\varphi Y^\uparrow)\varphi Z^\uparrow
+\left\{g\big(X^\uparrow,\big(\nablaM_{Y^\uparrow}\varphi\big)Z^\uparrow\big)-
         g\big(Y^\uparrow,\big(\nablaM_{X^\uparrow}\varphi\big)Z^\uparrow\big)\right\}\xi
\end{array}
$$
and hence
$$\begin{array}{rl}
R^M(W^\uparrow,Z^\uparrow,X^\uparrow,Y^\uparrow) & =R^N(W,Z,X,Y)\circ\pi
     -2g(X^\uparrow,\varphi Y^\uparrow)g(W^\uparrow,\varphi Z^\uparrow) +\\[1mm]
   & + g(Y^\uparrow,\varphi Z^\uparrow)g(W^\uparrow,\varphi X^\uparrow)-
     g(X^\uparrow,\varphi Z^\uparrow)g(W^\uparrow,\varphi Y^\uparrow).
 \end{array}
$$
Suppose that the base manifold $N$ satisfies Gray identities.
What are the corresponding curvature identities for the
upstairs manifold $M$?\\[2mm]
If $N$ is $(K_1)$ then
$$
\begin{array}{c}
R^M(X^\uparrow,Y^\uparrow,\varphi Z^\uparrow,\varphi W^\uparrow)-
            R^M(X^\uparrow,Y^\uparrow,Z^\uparrow,W^\uparrow)=\qquad\quad \\[1mm]
\qquad\qquad =-g(Y^\uparrow,W^\uparrow)g(Z^\uparrow,X^\uparrow)-
           g(Y^\uparrow,\varphi W^\uparrow)g(Z^\uparrow,\varphi X^\uparrow)\\[1mm]
\qquad\qquad\qquad
+g(X^\uparrow,W^\uparrow)g(Z^\uparrow,Y^\uparrow)+
           g(X^\uparrow,\varphi W^\uparrow)g(Z^\uparrow,\varphi Y^\uparrow).
\end{array}
$$
If $N$ is $(K_2)$ then
$$
\begin{array}{c}
R^M(\varphi X^\uparrow,Y^\uparrow,Z^\uparrow,W^\uparrow)+
R^M(X^\uparrow,\varphi Y^\uparrow,Z^\uparrow,W^\uparrow)+\qquad\qquad\\
\qquad\qquad +R^M(X^\uparrow,Y^\uparrow,\varphi Z^\uparrow,W^\uparrow)+
R^M(X^\uparrow,Y^\uparrow,Z^\uparrow,\varphi W^\uparrow)=0.
\end{array}
$$
If $N$ is $(K_3)$ then
$$
R^M(\varphi X^\uparrow,\varphi Y^\uparrow,\varphi Z^\uparrow,\varphi W^\uparrow)-
            R^M(X^\uparrow,Y^\uparrow,Z^\uparrow,W^\uparrow)=0.
$$
These relations are exactly the defined Gray identities for almost contact metric
manifolds for vector fields orthogonal to $\xi$.

\section{Properties and examples}
In their paper \cite{kn:JV81}, D. Janssens and L. Vanhecke have studied curvature
tensors for almost contact metric structures and defined {\em almost $C(\alpha)-$manifolds},
namely those almost contact metric manifolds whose curvature tensor
satisfies the following property:
$$
    \exists\alpha\in{\mathbf{R}}\ {\rm such\ that\ for\ all\ }X,Y,Z,W\in\chi(M)
$$
$$
   R(X,Y,Z,W)=R(X,Y,\varphi Z,\varphi W)+\alpha\left\{-g(X,Z)g(Y,W)+g(X,W)g(Y,Z)\right.
$$
$$
   +\left.g(X,\varphi Z)g(Y,\varphi W)-g(X,\varphi W)g(Y,\varphi Z)\right\}.
$$
This means that manifolds satisfying the first Gray identity $(K_{1\varphi})$ in the sense of Bonome et al.
are in fact $C(0)$-manifolds, while that manifolds satisfying $(G1)$ are
$C(1)$-manifolds. Note that cosymplectic, Sasakian and Kenmotsu manifolds are respectively
$C(0)$, $C(1)$ and $C(-1)$ manifolds (see Theorem 2.3, in \cite{kn:JV81}).

\medskip

Let us come back to Gray identities for an almost Hermitian manifold.
\medskip

It is known that
$K1\Rightarrow K2\Rightarrow K3$ (see \cite{kn:Gra76}, \S5).
Consequently we have
\begin{prop}
For a class ${\mathcal{L}}$ of almost contact metric manifolds,
denote by ${\mathcal{L}}_i$ the subclass of manifolds whose curvature
satisfies $Gi$, $i=1,2,3$. Then we have the following inclusions
$$
{\mathcal{L}}_1\subseteq{\mathcal{L}}_2\subseteq{\mathcal{L}}_3\subseteq{\mathcal{L}}.
$$
\end{prop}
As Gray remarked for Kaehlerian manifolds, we can say that {\em as i decreases, a manifold in
${\mathcal{L}}_i$ resembles Sasakian manifold more closely.}

\begin{prop}
Let $(M, \varphi,\xi,\eta,g)$ be a K-contact manifold satisfying G1 curvature identity.
Then the manifold $M$ is Sasakian.
\end{prop}
\proof
By using Proposition 7.5 in \cite{kn:Bla02}, p.94, a K-contact manifold whose curvature
satisfies $R_{XY}\xi=\eta(Y)X-\eta(X)Y$ is Sasakian. But this last relation is a consequence of
G1 identity. See also Proposition \ref{prop:cont_sas}.
\gata
\begin{prop}
Let $M$ be a contact metric manifold for which $\xi$ belongs to the
$(\kappa,\mu)$-nullity distribution, namely its curvature satisfies
\begin{equation}
\label{eq:kappa_mu}
R_{XY}\xi=\kappa\left(\eta(Y)X-\eta(X)Y\right)+\mu\left(\eta(Y)hX-\eta(X)hY\right)
\end{equation}
where $h=\frac12\ {\mathcal{L}}_\xi\varphi$ and $\kappa$, $\mu$ are constants.
Suppose $M$ satisfies $(G1)$ identity. Then $M$ is Sasakian.
\end{prop}
\proof
If $M$ is $(G1)$ then $R_{XY}\xi=\eta(Y)X-\eta(X)Y$ for all $X,Y\in\chi(M)$.
Combining with the fact that $\xi$ belongs to the $(\kappa,\mu)$-nullity distribution
we obtain
$$
  (\kappa-1)(\eta(Y)X-\eta(X)Y)+\mu(\eta(Y)hX-\eta(X)hY)=0
$$
for all $X,Y\in\chi(M)$. If $\mu\neq0$ this implies
$hY=\frac{1-\kappa}\mu\ Y$ for all $Y\in\ker\eta$.
We know that $h$ anticommutes with $\varphi$ and hence
one gets $\kappa=1$. But using Theorem 7.7, p. 103 in \cite{kn:Bla02}
it follows that $M$ is a Sasakian manifold.
If $\mu=0$ we immediately have $\kappa=1$.
\gata
\begin{prop}
Let $(M,\varphi,\eta,\xi,g)$ be a contact metric manifold satisfying $(G3)$ identity.
Then $M$ is K-contact.
\end{prop}
\proof
Choose a $\varphi -$adapted local orthonormal frame on $M$, namely
$\{X_i,\varphi X_i,\xi\}$, $i=1,\ldots, n$. Since $M$ is $(G3)$ the relation
$R(X,\xi,Y,\xi)=g(X,Y)$ holds for all $X,Y\in\ker\eta$. Taking $X=Y=X_i$
(respectively $X=Y=\varphi X_i$) one immediately obtains
$Ric(\xi,\xi)=2n$, where $Ric$ is the Ricci tensor on $M$.
Now we use the fact that {\em a contact metric manifold is K-contact if and only
if the Ricci tensor in the direction of the characteristic vector field $\xi$ is
equal to $2n$} (Theorem, p.65, \cite{kn:Bla76}).
\gata

\subsection{An example of almost contact metric manifold satisfying $G2$ but not $G1$.}
The generalized Heisenberg group $H(p,1)$ is defined as the set of matrices of real numbers
having the form
$$
a=\left[\begin{array}{lcr}
     1 & A & c\\
     0 & I_p & {}^tB\\
     0 & 0 & 1
  \end{array}\right]
$$
where $I_p$ is the identity $p\times p$ matrix, $A=(a_1,\ldots,a_p)$,
$B=(b_1,\ldots,b_p)\in{\mathbf{R}}^p$ and $c\in{\mathbf{R}}$. (Cf. \cite{kn:GC89}.)
$H(p,1)$ is connected, simply connected nilpotent Lie group of dimension $2n+1$.
We will consider $p=2$. A global system of coordinates $(x^1,x^2,y^1,y^2,z)$
on $H(2,1)$ is defined by $x^i(a)=a_i$, $y^i(a)=b_i$ for $i-1,2$ and $z(a)=c$.
The global vector fields defined by
$$
  X_i =2 \frac{\partial }{\partial x^i},\
  Y_i =2 \left(\frac{\partial }{\partial y^i} +x^i \frac{\partial}{\partial z} \right)\ {\rm for\ } i=1,2,\
  {\rm and\ }\xi = 2 \frac{\partial}{\partial z}$$
are left invariant. We take $\eta =\frac{1}{2} (dz-x^1 dy^1 - x^2 dy^2)$ and the metric
$$
  g=\frac{1}{4} (dx^1 \otimes dx^1 +dx^2 \otimes dx^2+ dy^1 \otimes dy^1 + dy^2 \otimes dy^2)+ \eta \otimes \eta.
$$
By direct computations we obtain that
$d \eta=-\frac{1}{2} (dx^1 \wedge dy^1 + dx^2 \wedge dy^2)$
and $\xi$ is the characteristic vector field, namely $\eta(\xi) =1$ and
$i_\xi d\eta =0$.
Moreover, the basis defined above is orthonormal:
$g(X_i, X_j)=g(Y_i,Y_j) = \delta _{ij}$,
$g(\xi, \xi) =1$ and $g(X_i, Y_j)=g(X_i, \xi)=g(Y_i, \xi)=0$.
One has $[X_i, X_j]=2\xi$ and the other brackets are equal to zero and therefore
it is easy to verify the Levi-Civita connection is given by the following formulas:
$$\nabla_\xi X_i = - Y_i = \nabla_{X_i} \xi$$
$$\nabla_\xi Y_i = X_i = \nabla_{Y_i} \xi$$
$$\nabla_{X_i} Y_i  = -\nabla_{Y_i} {X_i} = \xi $$
for $i=1,2$, the other derivatives being zero.
We compute also the Riemann-Christoffel curvature tensor field:
$$R(X_1, X_2, Y_1, Y_2) =-1,\qquad
  R(X_1, Y_2, X_2, Y_1) =-1$$
$$R(X_1, Y_1, X_2, Y_2) =-2\qquad
  R(X_i, Y_i, X_i, Y_i) =-3$$
$$R(X_i, \xi, X_i, \xi) =1\qquad
  R(Y_i, \xi, Y_i, \xi) =1 \: \: \: \: for \: \: \: i=1,2.$$
The other values are zero or can be obtained from these ones.
Define $\varphi$ by:
$$\varphi X_1 = \cos\theta Y_1 + \sin \theta Y_2\quad
  \varphi X_2= \sin \theta Y_1 - \cos \theta Y_2$$
$$\varphi Y_1 = -\cos \theta X_1 - \sin \theta X_2\quad
  \varphi Y_2 = -\sin \theta X_1 + \cos \theta X_2$$
and hence $(M,g, \varphi, \xi, \eta)$ is an almost contact metric manifold.
\begin{prop}
The structure is K-contact but not Sasakian.
\end{prop}
\proof For every $X,Y \in \chi (M)$ we have
$$g(\nabla_X \xi, Y) + g(\nabla_Y \xi, X)=0$$
On a Sasakian manifold, we should have $\nabla_X \xi = - \varphi X$
which implies in our case $\theta =0$. So, in general, $H(2,1)$ is not a Sasakian manifold.
\gata
\begin{prop}
On $H(2,1)$ the $G2$ identity holds, while $G1$ identity doesn't.
\end{prop}
Straightforward computations to prove $G2$. Moreover,
a K-contact manifold on which $G1$ holds is necessarily Sasakian.
This is not the case.
\gata

\subsection{Other examples.}
Let $(N,\bar g, J)$ be an almost Hermitian manifold. Consider the
warped product manifold $M={\mathbf{R}}\times_f N$, where $f=f(\theta)$
is the warping function and $\theta$ is the global parameter on ${\mathbf{R}}$.
Denote by $g=d\theta^2+f^2(\theta)\bar g$ the Riemannian metric on $M$.
Define also the global vector field $\xi=\frac\partial{\partial\theta}$,
the 1-form $\eta=d\theta$ and the $(1,1)$ tensor field $\varphi$ by
$\varphi X=JX$ if $X$ is tangent to $N$ and $\varphi\frac\partial{\partial\theta}=0$.
Thus $(\varphi,\xi,\eta,g)$ is an almost contact metric structure on $M$.
If $\bar\nabla$ and $\nabla$ are the Levi Civita connections on $N$, respectively
on $M$, we have
$$
\nabla_\xi X=\nabla_X\xi=\frac{f'}f\ X,\quad \nabla_\xi\xi=0,\quad
\nabla_XY=\bar\nabla_XY-ff'\bar g(X,Y)\xi,
$$
for all $X,Y$ tangent to $N$.\\[2mm]
The Riemann Christoffel curvature tensor is given by
\begin{equation}
\begin{array}{l}
R(W,\xi,X,Y)=0,\quad  R(W,\xi,X,\xi)=-\frac{f''}f\ g(X,W)\\[2mm]
R(W,Z,X,Y)=f^2\big[\bar R(W,Z,X,Y)+ \\[1mm]
 \qquad\qquad\qquad\qquad  +(f')^2\left(\bar g(X,Z)\bar g(Y,W)-\bar g(Y,Z)\bar g(X,W)\right)\big]
\end{array}
\end{equation}
In order to have one of the three curvature identities we immediately have
$$\frac{f''}f=-1$$
which implies that $f=\alpha\cos\theta+\beta\sin\theta$ with $\alpha$ and $\beta$ real constants.
At this one can state the following
\begin{prop}
The manifold $M$ is $G2$ (respectively $G3$) if and only if the almost Hermitian
manifold $N$ is $K2$ (respectively $K3$).
\end{prop}
\proof
One has the following relations:
$$
\begin{array}{l}
R(\varphi W,Z,X,\varphi Y)+R(W,\varphi Z,X,\varphi Y)+R(W,Z,\varphi X,\varphi Y)=\\[1mm]
\qquad  = f^2\big[ \bar R(JW,Z,X,JY)+\bar R(W,JZ,X,JY)+\bar R(W,Z,JX,JY) \big]\\[1mm]
\qquad \qquad  +(f')^2f^2\big( \bar g(X,Z)\bar g(Y,W)-\bar g(Y,Z)\bar g(X,W) \big)
\end{array}
$$
and
$$
R(W,Z,X,Y)-R(\varphi W, \varphi Z, \varphi X, \varphi Y)=
f^2\big[ \bar R(W,Z,X,Y)-\bar R(JW,JZ,JX,JY) \big].
$$
Hence the statement.
\gata

\begin{rem}
If $\dim N\geq 4$ then the manifold $M$ cannot be $G1$.
\end{rem}
\proof
Suppose $M$ satisfies $G1$ identity.
A straightforward computation gives
$$
\begin{array}{l}
\bar R(W,Z,JX,JY)-\bar R(W,Z,X,Y)=\left(1+(f')^2\right)\big[ \bar g(JX,W)\bar g(JY,Z)-\qquad\\[1mm]
\qquad\qquad -\bar g(JX,Z)\bar g(JY,W)+\bar g(Y,W)\bar g(X,Z)-\bar g(Y,Z)\bar g(X,W)\big].
\end{array}
$$
Since $f$ depends on $\theta$ (and it is not linear) while $\bar g$ and $\bar R$ do not,
it follows that $N$ is $K1$ and
$$
\bar g(JX,W)\bar g(JY,Z)-\bar g(JX,Z)\bar g(JY,W)+\bar g(Y,W)\bar g(X,Z)-\bar g(Y,Z)\bar g(X,W)=0
$$
for all $X,Y,Z,W$ tangent to $N$. This yields
\begin{equation}
\label{eq_for_G1}
\bar g(JY,Z)JX-\bar g(JX,Z)JY+\bar g(X,Z)Y-\bar g(Y,Z)X=0.
\end{equation}
If $\dim N\geq4$ we can choose $X$ and $Y$ such that $X$, $Y$,
$JX$ and $JY$ are linearly independent, so, the previous equality is impossible.
\gata
\begin{ex}\rm
On $M={\mathbf{R}}^4\times (-\pi/2,\pi/2)$ consider the global coordinates
$x$, $y$, $u$, $v$ and $z$. Consider the Riemannian metric
$g=dz^2+\cos^2 z\left(dx^2+dy^2+du^2+dv^2\right)$ and the almost contact
structure defined by: $\xi=\partial_z$, $\eta=dz$, $\varphi\partial_u=\partial_v$,
$\varphi\partial_v=-\partial_u$ and $\varphi\partial_z=0$.
Then $M$ is $G2$ but not $G1$.\\[2mm]
Similarly for $M={\mathbf{R}}^4\times (0,\pi)$ and
$g=dz^2+\sin^2 z\left(dx^2+dy^2+du^2+dv^2\right)$.

\medskip
This kind of structure is called {\em sine-cone} and gives way to construct many
geometric objects (e.g. nearly Kaehler structures starting from a 5-dimensional
Sasaki Einstein manifold). Cf. \cite{kn:FIMU06}.
\end{ex}
\begin{prop}
Let $N$ be a surface (which is automatically Kaehler)
and consider the warped product manifold $M$ as above.
Then $M$ satisfies $G1$.
\end{prop}
\proof
The statement follows from the fact that a Kaehler manifold is $K1$
and the equation (\ref{eq_for_G1}) is satisfied in dimension 2.
\gata

\subsection{Hypersurfaces of almost Hermitian manifolds}
Let $(\widetilde M, J, \widetilde g)$ a $(2n+2)$-dimensional Kaehler manifold,
and let $M$ be a totally umbilical (real) hypersurface in $\widetilde M$.
Denoting by $N$ the unit normal on $M$ and let $A$ be the Weingarten
operator and $h$ the scalar second fundamental form.
As $M$ is totally umbilical, we have that $AX = \beta X$, for all $X$ tangent to $M$,
with $\beta \in C^\infty (M)$.
\medskip

It is well known the fact that on $M$ we can define an almost contact metric structure
(see e.g. \cite{kn:Bla02}). More precisely,
we take $\xi =-JN$ and for $X\in \chi(M)$ we decompose $JX$ as:
$$JX=\varphi X + \eta(X)N.$$
Let $g$ be the restriction of the metric $\widetilde g$ on $M$.
Denote by $\widetilde \nabla$ (respectively $\nabla$) the Levi-Civita connection on $\widetilde M$
(respectively on $M$). Then, by the formula of Gauss,
%  $$\widetilde \nabla_X Y =\nabla_X Y +h(X,Y)N$$
% In particular
one has
$$\widetilde \nabla_X \xi = \nabla_X \xi +h(X,\xi)N$$
On the other hand we have $\widetilde \nabla_X \xi= -J \widetilde \nabla_X N=JAX=\varphi AX+\eta(AX)N$.
Hence
$$
  \nabla_X \xi =\varphi AX\quad {\rm and}\quad
  h(X,\xi)=\eta(AX).
$$
Suppose now that $M$ satisfies the $(G3)$ identity.
This implies that
% $$R(X,\xi, Y,\xi) =g(X,Y)-\eta(X)\eta(Y) \:\:\:\:\: \forall X,Y \in \chi(M).$$
\begin{equation}
\label{ec0}
R(X,\xi, Y,\xi) =g(X,Y) \:\:\:\:\: \forall X,Y \in Ker(\eta)
\end{equation}
We should compute $R(X,\xi)\xi = \nabla_X\nabla_\xi \xi-\nabla_\xi \nabla_X \xi-\nabla_{[X,\xi]} \xi$.
Since $M$ is totally umbilical, we have that $\nabla_X \xi=\beta \varphi X$.
Thus $\nabla_\xi \xi =0$.
Then
$$\nabla_\xi \nabla_X \xi = \xi(\beta) \varphi X + \beta(\nabla_\xi \varphi) X + \beta \varphi \nabla_\xi X.$$
But
$$\nabla_\xi X = \beta \varphi X -[X,\xi]$$
and so
$$\nabla_\xi \nabla_X \xi = \xi(\beta) \varphi X + \beta(\nabla_\xi \varphi) X +
           \beta^2  \varphi^2 X -\beta \varphi [X,\xi].
$$
It follows that
$$R(X,\xi)\xi = -\xi(\beta) \varphi X - \beta (\nabla_\xi \varphi) X + \beta ^2 X.$$

Now, due the fact $M$ is Kaehler, we have that
$$\widetilde \nabla (JY) = J \widetilde \nabla_X Y=J(\nabla_X Y + h(X,Y)N)=\varphi \nabla_X Y + \eta(\nabla_X Y)N-h(X,Y)\xi$$
On the other hand
$$
   \widetilde \nabla (JY)=\widetilde \nabla_X (\varphi Y + \eta(Y)N)=
               \nabla(\varphi Y) +h(X, \varphi Y)N + X \eta(Y)N-\eta(Y)\beta X.
$$
Identifying the tangent and the normal parts of $\widetilde \nabla (JY)$ we obtain respectively
\begin{equation}
\label{ec1}
(\nabla_X \varphi )Y= \beta \eta(Y)X - \beta g(X,Y)\xi
\end{equation}
\begin{equation}
(\nabla_X \eta)(Y)=-\beta g(X,\varphi Y).
\end{equation}
Putting $X =\xi$ in (\ref{ec1}) we have
$(\nabla_\xi \varphi)Y = \beta \eta(Y) \xi -\beta g(\xi, Y)\xi =0$
which implies
$$\nabla_\xi \varphi =0.$$
Then
$$
   R(X,\xi) \xi= -\xi(\beta) \varphi X + \beta ^2 X
$$
From (\ref{ec0}) we have that
$$g(\beta^2 X - X -\xi(\beta) \varphi X , Y) =0,\quad \forall Y \in \ker(\eta).$$
As $X$ and $\varphi X$ are linearly independent (and belong to $\ker (\eta)$),
we obtain that $\beta =\pm 1$.

We obtain that
$$
    AX =\pm X, \:\:\:\: and \:\:\:\: h(X,Y) = \pm g(X,Y).
$$
For $\beta =-1$ it follows that
$$(\nabla_X \varphi)Y = g(X,Y) \xi -\eta(Y)X.$$
By Theorem 6.14 in \cite{kn:Bla02} this implies that $M$ is Sasakian.
\begin{prop}
Let $M$ a totally umbilical hypersurface of a Kaehler manifold $\widetilde M$ endowed with the
usual almost contact metric structure. If $M$ satisfies the $G3$ identity, then $M$ is a Sasakian manifold
and hence $M$ satisfies all $Gi$, for $i=1,2,3$.
\end{prop}

More generally, if the second fundamental form of $M$ is given by
$$
  h(X,Y)=\lambda\eta(X)\eta(Y)+\mu g(X,Y),\quad\forall X,Y\in\chi(M)
$$
with $\lambda$ and $\mu$ smooth functions on $M$, i.e. $M$ is {\em totally
quasi umbilical}, and if $M$ satisfies $(G3)$ identity, then it is Sasakian.
As consequence, there is no {\em cylindrical submanifold} satisfying $(G3)$ and
whose second fundamental form being $h(X,Y)=\lambda\eta(X)\eta(Y)$.

\medskip
{\bf Acknowledgements.}
The authors would like to thank Professor A. Bonome for discussions we had in Santiago de Compostela
concerning this subject. They also wish to express their gratitude to Professor D.E. Blair for
useful comments and suggestions during the preparation of this paper.

\vspace{3mm}

{\scriptsize
Authors' addresses:\\[2mm]
Raluca Mocanu,\\
University of Bucharest,
Faculty of Mathematics,
Str. Academiei n.14,
s 1, Bucharest,
Romania\\
e-mail: xipita@yahoo.com\\[2mm]
Marian Ioan Munteanu,\\
University 'Al.I.Cuza' of Ia\c si,
Faculty of Mathematics,
Bd. Carol I, no.11,
700506 Ia\c si,
Romania\\
e-mail: munteanu@uaic.ro
}

\end{document}